\documentclass[10pt,a4paper]{amsart}

\usepackage{amsmath}

\DeclareMathOperator{\dimH}{\mathrm{dim}_\mathrm{H}}

\begin{document}

\title{On a problem by R. Salem concerning Minkowski's question mark function}
\author{Tomas Persson}
\address{Tomas Persson, Centre for Mathematical Sciences, Lund University, Box 118, 22100 Lund, Swed
en}
\email{tomasp@maths.lth.se}
\urladdr{http://www.maths.lth.se/~tomasp}

\date{\today}

\begin{abstract}
Minkowski's question mark function is strictly increasing on $[0, 1]$
and hence defines a Stieltjes measure on $[0, 1]$. A problem
originating from Salem in 1943, is to determine whether the Fourier
series of this measure decay to zero or not. The purpose of this note
is to mention that a recent result by Jordan and Sahlsten implies that
the Fourier transform decays to zero with a polynomial speed.
\end{abstract}

\maketitle 

Minkowski's question mark function $?$ was introduced by H. Minkowski
in order to give an example of a one-to-one map between the rational
in $(0,1)$ and the quadratic irrationals in $(0,1)$. R. Salem showed
in \cite{salem} that
\[
? (x) = \sum_{k=1}^\infty \frac{(-1)^{k-1}}{2^{a_1(x) + a_2(x) +
    \cdots + a_k (x) - 1}},
\]
where $[0; a_1(x), a_2(x), \ldots]$ is the continued fraction
expansion of $x \in [0,1]$. The function $? \colon [0,1] \to [0,1]$ is
strictly increasing and so it defines a Stieltjes measure $\mu$ of
mass one.  It was mentioned in \cite{salem} that it does not seem to
be known whether the Fourier series of $\mu$ decays to zero at
infinity.

According to J. R. Kinney \cite{kinney} the Hausdorff dimension of the
measure $\mu$ is given by
\[
\dimH \mu = \biggl(2 \int_0^1 \log_2 (1+x) \, \mathrm{d} \mu (x)
\biggr)^{-1}.
\]
Since $0 < \log_2 (1+x) < 1$ on $(0,1)$ we clearly have $\dimH \mu >
1/2$.

The measure $\mu$ can also be obtained in the following equivalent
way, as mentioned by F.~Canto-Mart\'{i}n, H.~Hedenmalm, and
A.~Montes-Rodr\'{i}guez on page 43 of \cite{hedenmalm}. Let $G \colon
[0,1) \to [0,1)$ be the Gau\ss-map, defined by
\[
G \colon x \mapsto 1/x \mod 1.
\]
The map $G$ is continuous and decreasing on each of the intervals
\[
I_k = (1/(k+1), 1/k ], \qquad k = 1, 2, 3, \ldots
\]
The measure $\mu$ is then the $G$-invariant Bernoulli measure obtained
by defining
\[
\mu(I_k) = 2^{-k}.
\]

The problem concerning the Fourier series of $\mu$ mentioned by
R. Salem is also mentioned in the above mentioned paper by
F. Canto-Mart\'{i}n, H. Hedenmalm, and A. Montes-Rodr\'{i}guez.

Since $\dimH \mu > 1/2$, as mentioned above, it is an immediate
consequence of the result of T. Jordan and T. Sahlsten \cite[Theorem
  1.3]{jordansahlsten} that there is a number $\eta > 0$ such that
\[
|\hat{\mu} (\xi)| = O(|\xi|^{-\eta}) \qquad \text{as } |\xi| \to
\infty.
\]

Hence, the theorem by Jordan and Sahlsten solves the problem raised by
R. Salem.

I noticed this on January 1, 2015 and told T. Jordan and T. Sahlsten
about it the same day.  The day before, December 31, 2014, another
solution to Salem's problem was submitted to Arxiv by S. Yakubovich
\cite{yakubovich}, and appeared on Arxiv on January 5, 2015. Before
January 5, I did not know about Yakubovich's result. The paper by
T. Jordan and T. Sahlsten, from which the result immediately follows,
appeared more than one year before, on December 12, 2013.

\end{document}